\documentclass[twocolumn]{autart} 
\usepackage{amsmath,amssymb,euscript,yfonts,psfrag,latexsym,graphicx}
\usepackage{bbm,color,amstext,wasysym,parskip}

\graphicspath{{./},{./figures/}}
\usepackage{caption}
\usepackage{subcaption}
\usepackage{float}
\usepackage{array}
\usepackage{url}
\usepackage{algorithm}
\usepackage{algorithmic}
\usepackage{xcolor}

\usepackage{algorithmic}
\usepackage{caption}
\usepackage{algorithm}

\newcommand{\mE}{{\mathbb E}}

\newcommand{\mR}{{\mathbb R}}

\newcommand{\cL}{{\mathcal L}}

\newcommand{\cP}{{\mathcal P}}
\newcommand{\cQ}{{\mathcal Q}}

\newcommand{\cY}{{\mathcal Y}}

\newcommand{\Tr}{\operatorname{Tr}}

\newcommand{\argmin}{\operatorname{argmin}}
\usepackage{ulem}

\definecolor{grey}{rgb}{0.6,0.6,0.6}
\definecolor{lightgray}{rgb}{0.97,.99,0.99}
\newcommand{\comment}[1]{\textcolor{grey}{#1}}

\let\classAND\AND
\let\AND\relax
\usepackage{algorithmic}

\let\AND\classAND
\AtBeginEnvironment{algorithmic}{\let\AND\algoAND}

\floatname{algorithm}{Algorithm}

\begin{document}
\begin{frontmatter}
\title{An optimal control approach to particle filtering}

\thanks[footnoteinfo]{This work was supported by the NSF under grant 1942523 and 2008513.}

\author[qinsheng]{Qinsheng Zhang}\ead{qzhang419@gatech.edu},  
\author[amirhossein]{Amirhossein Taghvaei}\ead{amirtag@uw.edu}, 
\author[yongxin]{Yongxin Chen}\ead{yongchen@gatech.edu}  %

\address[qinsheng]{ Machine Learning Center, Georgia Institute of Technology, Atlanta, GA, USA}             %
\address[amirhossein]{Department of Aeronautics and Astronautics, University of Washington, Seattle, WA, USA}  
\address[yongxin]{School of Aerospace Engineering, Georgia Institute of Technology, Atlanta, GA, USA}        %


\begin{abstract}
We present a novel particle filtering framework for continuous-time dynamical systems with continuous-time measurements. Our approach is based on the duality between estimation and optimal control, which allows reformulating the estimation problem over a fixed time window into an optimal control problem. The resulting optimal control problem has a cost function that depends on the measurements and the closed-loop dynamics under optimal control coincides with the posterior distribution over the trajectories for the corresponding estimation problem. This type of stochastic optimal control problem can be solved using a remarkable technique known as path integral control.
By recursively solving these optimal control problems using path integral control as new measurements become available we obtain an optimal control based particle filtering algorithm. A distinguishing feature of the proposed method is that it uses the measurements over a finite-length time window instead of a single measurement for the estimation at each time step, resembling the batch methods of filtering, and improving fault tolerance. The efficacy of our algorithm is illustrated with several numerical examples.
\end{abstract}

\end{frontmatter}

\section{Introduction}
In control engineering, filtering refers to estimating the true state of a dynamical system using the the raw sensor measurements. It is a critical component in feedback control and plays an indispensable role in almost all applications related to control.  There are many theories and algorithms for filtering that have been developed. The celebrated Kalman filter is for linear dynamics driven by Gaussian noise. It is optimal in the sense of mean-square error. It also computes the exact posterior distribution of the state given the available measurements. For nonlinear systems, the filtering problem is much more challenging; the posterior distribution of the state rarely has a simple parametrization. To attain the posterior distribution, one needs to solve a stochastic partial differential equation known as the Kushner–Stratonovich equation. The methods relying on discretizing the state space and  the Kushner–Stratonovich equation are computationally infeasible for high dimensional problems. There are some algorithms that approximate posterior distributions using Gaussian distributions, including the extended Kalman filter (EKF) and the unscented Kalman filter (UKF). However, the performance of this type of methods deteriorates as the posterior distribution drifts away from the Gaussian family. 

One approach that avoids brutal force discretization of the state space while still retaining the richness of the posterior distributions is representing the distributions with particles. This type of methods are known as the particle filtering. Over the last decades, many different versions of particle filtering algorithms have been proposed~\cite{TagMehMey20}~\cite{GuaJohLee17}~\cite{RuiKap17}~\cite{SarSot08}. In the standard setup of particle filtering, the posterior distribution at the current step is approximated by $K$ weighted particles. These particles are propagated forward following a proposal density and then combined with the next measurement to estimate the posterior distribution at the next time step. The implementation of particle filtering is extremely easy if the proposal density is simple, which makes particle filtering a popular method for nonlinear filtering. Theoretically, it can be shown that as the number of particles $K$ goes to infinity, the empirical distribution of the particles converges to the true posterior distribution at each time step in some suitable sense. In practice, however, due to the potentially large difference between prior dynamics and posterior dynamics, the weights of the particles become degenerate quickly~\cite{DouBriSen06}. That is, the weights of most of the particles become negligible and the mass of the particles only concentrates on a few particles, rendering a small effective particle size. A resampling step is commonly adopted to mitigate the effects of degenerate weights. However, both in theory and in practice, how to choose a proper proposal density is critical and most particle filtering algorithms still perform poorly in high-dimensional problems \cite{BenBicLi08,BesCriJasWhi14}, largely due to particle degeneracy.
  
 In this work, we consider nonlinear filtering problems for continuous-time diffusion dynamics with continuous-time measurements. We present a new particle filtering method based on an elegant duality between estimation and optimal control \cite{MitNew03,KimMeh20}. Building on this duality, we are able to obtain a superior proposal density by (approximately) solving an optimal control problem and thus establish a particle filtering algorithm with great performance. Moreover, this duality makes it natural to resample the particles from the past; this is different to most particles filtering algorithms only samples in the present. This extra flexibility of updating samples in the past provides us the opportunity to correct numerical errors or errors induced by outlier in the previous filtering steps and makes the algorithm more robust to mistakes and outlier measurements. Empirically, we also observe that extending the sampling to the past, with proper proposals, can significantly mitigate the particle degeneracy issue.

The proposed algorithm is most related to those proposed in \cite{DouBriSen06}, \cite{RuiKap17}, and \cite{Bal09}. In \cite{DouBriSen06} a block sampling strategy is proposed to resample particles in the past as in our algorithm. However, they focus on an abstract framework for general discrete-time systems. How to leverage the structure of the underlying dynamics to construct a proper proposal distribution is not discussed explicitly. In \cite{RuiKap17}, an optimal control approach to smoothing is proposed. However, they consider smoothing problem over a fixed-time window. Moreover, though the dynamics they use is continuous-time diffusion, their measurement model is discrete-time. The same setting with continuous-time diffusion and discrete-time measurement is used in \cite{Bal09}. In addition, even though some path integral idea is used, the algorithm in \cite{Bal09} is grid-based, not particle based. There are also some other particle filtering algorithms such as feedback particle filtering \cite{YanMehMey13,TagMehMey20} and particle flow filter \cite{DauHuaNou10} that aim to improve the performance by using a better proposal.

The rest of the paper is structured as follows. In Section \ref{sec:background}, we provide a brief introduction to particle filtering and stochastic optimal control. The remarkable duality between filtering and optimal control is presented in Section \ref{sec:duality}. We then use this optimal control formulation of filtering to derive our particle filtering algorithm in Section \ref{sec:controlfiltering}. The algorithm is illustrated in Section \ref{sec:examples} with several numerical examples. This is followed by a concluding remark in Section \ref{sec:conclusion}. 

\section{Background}\label{sec:background}
In this section we provide some background knowledge on particle filtering as well as optimal control that is closely related to our proposed method.
\subsection{Particle filtering}
The standard setting of particle filtering is over a discrete-time dynamic system 
	\begin{subequations}
	\begin{equation}
		x_{t+1} \sim p_X(x_{t+1} \mid x_t), \quad x_0 \sim p_X(x_0)
	\end{equation}
with observation model
	\begin{equation}
		y_t \sim p_Y(y_t \mid x_t).
	\end{equation}
	\end{subequations}
Here $p_X(\cdot), p_X (\cdot | \cdot)$ and $p_Y(\cdot |\cdot)$ denote the prior distribution of the initial state, transition probability, and measurement probability respectively. 
Continuous-time systems can be converted into this form via a discretization over time. 

The main idea of particle filtering is to represent the posterior distributions $p(x_t | y_{1:t})$ by a collection of particles $\{x_t^1,x_t^2, \ldots,x_t^K\}$ and the corresponding weights $\{w_t^1,w_t^2,\ldots,w_t^K\}$ satisfying $\sum_{i=1}^K w_t^i=1$. More specifically, the posterior distribution is approximated by
    \begin{equation}
        p(x_t | y_{1:t}) \approx \sum_{i=1}^K w_t^i \delta_{x_t^i}
    \end{equation}
where $\delta_x$ is the Dirac delta distribution located at $x$.
One of the most widely used particle filtering algorithm is the sequential important resampling (SIR) particle filter \cite{DouBriSen06}. It starts with $K$ independent samples $\{x_0^i\}_{i=1}^K$ from the prior distribution $p_X(\cdot)$. Since these samples are independently sampled, they are assigned equal weights, that is, $w_0^1=w_0^2=\cdots=w_0^K=1/K$. SIR uses the following updates 
\begin{align*}
    x^i_{t+1} &\sim p_X(x_{t+1}|x^i_t)\\
    w^i_{t+1} &= \frac{w^i_t p_Y(y_{t+1}|x^i_{t+1})}{\sum_{j=1}^K w^j_t p_Y(y_{t+1}|x^j_{t+1})}. 
\end{align*}
to iteratively approximate $p(x_t|y_{1:t})$.
A resampling step is implemented after several steps to avoid the weight degeneracy. The weight degeneracy is quantified by the effective ratio
	\[
    		\gamma = \frac{1}{K \sum_{i=1}^K (w_t^i)^2}.
	\]
The value of $\gamma$ has maximum $1$, achieved when the weights are uniform. When the effective ratio is below a certain threshold $\gamma_{thres}$, a resampling step is carried out. 
During resampling, $K$ independent new samples are generated from the weighted discrete distribution $\tilde{x}^i_t \sim \sum_{i=1}^K w^i_t \delta_{x^i_t}$. After resampling, the posterior distribution is approximated with empirical distribution of the new samples with equal weights, i.e., $p(x_t|y_{1:t}) \approx \frac{1}{K}\sum_{i=1}^K \delta_{\tilde{x}^i_t}$.

\subsection{Stochastic optimal control}
Consider the stochastic dynamics described by the stochastic differential equation (SDE) \cite{SarSol19}
	\begin{equation}\label{eq:sysSDE}
		dX_t = b(t,X_t) dt + \sigma(t,X_t) (u_tdt+dW_t)
	\end{equation}
where $X_t\in \mR^n, u_t\in \mR^m$ denotes the state and control input respectively, and $W_t\in\mR^m$ represents a standard Wiener process. The drift $b(\cdot,\cdot)$ and the input channel matrix $\sigma(\cdot,\cdot)$ are assumed to be Lipschitz continuous and bounded. 

In the finite horizon stochastic optimal control problem \cite{LewVraSyr12} one seeks an optimal feedback control strategy that minimizes the cost function
	\begin{equation}\label{eq:cost}
		J(u) = \mE \left\{\int_{0}^T w(t, X_t, u_t) dt + \Psi(X_T)\right\}
	\end{equation}
over a fixed time interval $[0,\, T]$. Here, $w$ and $\Psi$ represent running cost and terminal cost respectively. 
This problem can be solved via dynamic programming \cite{Ber95,Eva98}, which boils down to solving the Hamilton-Jacobi-Bellman (HJB) equation \cite{Eva98}
	\begin{equation}\label{eq:HJB}
		\frac{\partial V_t}{\partial t} + \min_{\alpha\in \mR^m} \{\cL_t^\alpha V_t +w(t,x,\alpha)\} =0,\quad V_T(\cdot) = \Psi(\cdot),
	\end{equation}
where $\cL_t^\alpha$ denotes the generator of the controlled process \eqref{eq:sysSDE} defined as
	\begin{align}
		\cL_t^\alpha f(x) = (b(t,x)+\sigma(t,x)\alpha)\cdot \nabla f(x) +\\ \frac{1}{2}\Tr(\sigma(t,x)\sigma(t,x)'\nabla^2 f(x)) \nonumber
	\end{align}
for any sufficiently smooth $f(\cdot)$. The space-time function $V_t(x)$ is known as the cost-to-go function, capturing the minimum cost over the time window $[t,\,T]$ conditioned on $X_t = x$.
The optimal control strategy is of state feedback form $u_t^\star = \alpha^\star(t,X_t)$ with
	\begin{equation}\label{eq:Lta}
		\alpha^\star(t,x) = \argmin_{\alpha\in \mR^m} \{\cL_t^\alpha V_t +w(t,x,\alpha)\}.
	\end{equation}

The filtering algorithm developed in this work is closely related to the special case of stochastic control problems where the running cost is of the form
	\begin{equation}\label{eq:costspecial}
		w(t,x,\alpha) = g(t,x)+\frac{1}{2}\|\alpha\|^2.
	\end{equation}
Clearly, with this running cost, the minimization in \eqref{eq:HJB} can be solved in closed-form, yielding the optimal policy
	\begin{equation}\label{eq:optcontrol}
		u_t^\star = -\sigma(t,X_t)'\nabla V_t(X_t),
	\end{equation}
and the HJB equation \eqref{eq:HJB} simplifies to 
	\begin{equation}\label{eq:HJBquadratic}
		\frac{\partial V_t}{\partial t} + b\cdot \nabla V_t+g-\frac{1}{2}\nabla V_t'\sigma\sigma' \nabla V_t+ \frac{1}{2}\Tr(\sigma\sigma'\nabla^2V_t)\!=\!0.
	\end{equation}

The running cost \eqref{eq:costspecial} plays a crucial role in our framework. The quadratic cost in control in \eqref{eq:costspecial} quantifies the difference of the controlled and the uncontrolled ($u_t\equiv 0$) process. More specifically, denote $\cP_u$ the measure over the path space $\Omega = C([0,T]; \mR^n)$ induced by the dynamics \eqref{eq:sysSDE}, and $\cP_0$ the measure associated with the uncontrolled process, then by the celebrated Girsanov theorem \cite{SarSot08}, 
	\begin{equation}\label{eq:girsanovthm}
		\frac{d\cP_u}{d\cP_0} = \exp \left\{\int_0^T \frac{1}{2}\|u_t\|^2 dt + u_t'dW_t\right\}.
	\end{equation}
It follows that the Kullback-Leibler divergence between $\cP_u$ and $\cP_0$ is \cite{SarSot08}
	\begin{equation}\label{eq:KL}
		{\rm KL} (\cP_u\,\|\, \cP_0):= \int_\Omega d\cP_u \log\frac{d\cP_u}{d\cP_0} = \mE \left\{\int_0^T\frac{1}{2}\|u_t\|^2 dt\right\},
	\end{equation}
where the expectation is with respect to the controlled process. Thus, the optimal control problem with running cost \eqref{eq:costspecial} can be equivalently written as
	\begin{equation}\label{eq:cPuform}
		\min_{\cP_u} \mE_{\cP_u} \left\{\int_0^T g(t,X_t)dt + \Psi(X_T)\right\} + {\rm KL}(\cP_u\,\|\,\cP_0).
	\end{equation}
Note that the optimization variable becomes $\cP_u$ instead of the control policy; the two are equivalent as the control policy fully determines the measure $\cP_u$ and vice versa \cite{ThiKap15}.

\subsection{A linear approach to stochastic optimal control}\label{sec:linearapproach}
When the cost is of the form \eqref{eq:costspecial}, it turns out that the above nonlinear optimal control problem can be solved in a linear manner \cite{KapRui16,ThiKap15,WilAldThe17,SarSot08,DouBriSen06,ZhaWanSch14,ThaKapTotGom20,BerHenDouJac19,Rei18,GuaJohLee17,RuiKap17,HenBisDou17}.  One way to see it is through the logarithmic transformation \cite{FleRis75} of the HJB equation \eqref{eq:HJBquadratic}. More specifically, let 
	\begin{equation}
		V_t(x) = -\log \varphi(t,x), 
	\end{equation}
then a straightforward calculation points to
	\begin{equation}\label{eq:varphipde}
		\frac{\partial \varphi}{\partial t} + b\cdot \nabla \varphi - g \varphi +\frac{1}{2}\Tr (\sigma\sigma' \nabla^2 \varphi) = 0,~ \varphi(T,\cdot) = \exp\{-\Psi\}.
	\end{equation}
The associated optimal control strategy reads
	\begin{equation}
		u_t^\star = \sigma(t,X_t)'\nabla\log\varphi(t,X_t).
	\end{equation}
Note that unlike the HJB \eqref{eq:HJBquadratic} which is nonlinear, \eqref{eq:varphi} is a linear partial differential equation (PDE); it is the Backward Kolmogorov equation \cite{SarSol19} associated with the (uncontrolled $u_t\equiv 0$) process \eqref{eq:sysSDE} and killing rate $g$. This transformation is remarkable; linear PDE is often much easier to solve than nonlinear PDE. In fact, this special PDE \eqref{eq:varphipde} can be solved through Monte Carlo sampling as discussed below. This idea is the foundation of the path integral control \cite{ThiKap15}.

To distinguish the processes associated with different control policies, we denote by $X^u$ the diffusion process \eqref{eq:sysSDE} with feedback control $u$ and $\mE_{X^u}$ the corresponding expectation. By the celebrated Feynman-Kac formula \cite{SarSol19}, we have 
	\begin{subequations}\label{eq:varphi}
	\begin{equation}
		\varphi(t,x) \!=\! \mE_{X^0}\! \left\{\!\exp\{-\int_t^T g(t,X_\tau^0)d\tau - \Psi(X_T^0)\}\mid X_t^0 = x\right\}.
	\end{equation}
Moreover, the optimal control at $(t,x)$ is \cite{ThiKap15}
	\begin{align}
		&u^\star (t,x) = \sigma(t,x)'\nabla \log \varphi(t,x)=\\ 
        &\lim_{s\searrow t} \frac{\mE_{X^0}\{\exp\{-\int_t^T g(\tau,X_\tau^0)d\tau - \Psi(X_T^0)\}\int_t^s dW_\tau\!\mid\! X_t^0 = x\}}{(s-t)\mE_{X^0}\{\exp\{-\int_t^T g(\tau,X_\tau^0)d\tau - \Psi(X_T^0)\}\!\mid\! X_t^0 = x\}}. \nonumber
	\end{align}
	\end{subequations}
We remark that the expectation is with respect to the uncontrolled process $X^0$. This is counter-intuitive; it says that one can recover the optimal control $u^\star$ by taking expectation with respect to the process with zero control. In practice, one can simulate many trajectories from $X^0$ with $X_t^0 = x$ and approximate $\varphi$ and $u^\star$ by taking average over these trajectories. This is exactly the path integral control \cite{ThiKap15}.

One potential issue of \eqref{eq:varphi} is that the variance of the Monte Carlo approximation could be high and thus a large number of samples is needed to achieve a reasonable accuracy. This drawback can be mitigated by importance sampling as follows \cite{ThiKap15}. Let
	\begin{subequations}\label{eq:Su}
	\begin{equation}
		S^u(t) = \int_t^T [g(\tau,X_\tau^u)+\frac{1}{2}\|u_\tau\|^2]d\tau+\int_t^T u_\tau'dW_\tau + \Psi(X_T^u),
	\end{equation}
then by the Girsanov theorem \eqref{eq:girsanovthm} we obtain
	\begin{equation}\label{eq:girsanov}
		\exp[-S^u(t)] = \exp\{-\int_t^T g(\tau,X_\tau^u)d\tau - \Psi(X_T^u)\}\frac{d\cP_0}{d\cP_u}.
	\end{equation}
	\end{subequations}
Plugging \eqref{eq:Su} into \eqref{eq:varphi} yields
	\begin{subequations}\label{eq:important}
	\begin{equation}\label{eq:important1}
		\varphi(t,x) = \mE_{X^u} \left\{\exp\{-S^u(t)\}\mid X_t^u = x\right\},
	\end{equation}
and
	\begin{eqnarray}\label{eq:important2}
		u^\star (t,x) -u(t,x) =&&\\
		&&\hspace{-1cm}\lim_{s\searrow t}\! \frac{\mE_{X^u}\{\exp\{-S^u(t)\}\int_t^s dW_\tau\!\mid\! X_t^u = x\}}{(s-t)\mE_{X^u}\{\exp\{-S^u(t)\}\!\mid\! X_t^u = x\}}.\nonumber
	\end{eqnarray}
	\end{subequations}
In practice, one samples $N$ trajectory starting from $X^u_{t,i}=x$ to calculate
\begin{align*}
    w_i&= \exp\left\{\!\!-\!\!\int_t^T \![g(\tau,X_{\tau,i}^u)\!+\! \frac{1}{2}\|u_{\tau,i}\|^2]d\tau\!\right.
    \\&\left.-\!\int_t^T \!\!u_{\tau,i}'dW_{\tau,i} \!-\! \Psi(X_{T,i}^u)\!\right\}.
\end{align*}
Here $W_{t,i}, u_{t,i}$ denote the noise and control signal associated with the trajectory $X^u_{t,i}$, respectively. 
The estimated optimal control then becomes
	\begin{equation}\label{eq:estimation}
		\hat u^\star(t,x) = u(t,x) + \sum_{i=1}^N\frac{w_i \Delta W_{t,i}}{\Delta t \sum w_j}
	\end{equation}
where $\Delta t$ is the stepsize and $\Delta W_{t,i}$ is one discretized step of $W_{t,i}$.
The variance of estimation \eqref{eq:estimation} is governed by the optimality of the policy $u(t,x)$. When $u$ is close to the optimal policy $u^\star$, the sample variance is small \cite{ThiKap15}. In fact, when $u$ coincides with the optimal control $u^\star$, the sample variance is $0$, meaning one can estimate \eqref{eq:estimation} with one trajectory.

\section{Smoothing as stochastic control}\label{sec:duality}
Consider a diffusion process with measurement noise
	\begin{subequations}\label{eq:filteringdyn}
	\begin{align}\label{eq:filteringdyn1}
		dX_t &= b(t,X_t) dt + \sigma(t,X_t) dW_t,\quad X_0\sim \nu_0
		\\\label{eq:filteringdyn2}
		dY_t &= h(t,X_t)dt + \sigma_B dB_t,\quad Y_0 = 0
	\end{align}
	\end{subequations}
where the measurement $Y_t\in \mR^p$ is corrupted by white noise $dB_t$ weighted by $\sigma_B>0$ and the initial state $X_0$ follows some prior distribution $\nu_0$. 
The smoothing problem is a particular type of Bayesian inference problem that aims at estimating the distribution of $X_t$ for $0\le t \le T$ given the full history of measurement $\{Y_t,\, 0\le t\le T\}$.  

It was discovered in \cite{MitNew03,KimMeh20} that the smoothing problem can be reformulated as a stochastic optimal control problem whose cost function depends on the measurements. To see this, denote the measure over the path space $\Omega$ induced by the process \eqref{eq:filteringdyn1} by $\cP$. This serves as the prior measure for this Bayesian inference problem. Denote the posterior distribution over $\Omega$ by $\cQ^Y$. 
By Kallianpur-Striebel formula \cite{Kle05},
\begin{align}\label{eq:KS}
    \frac{d\cQ^Y}{d\cP}\!\propto\!
    \exp\!\left\{
       \! -\frac{1}{\sigma_B^2}
        \left[
            \int_0^T [Y_t dh+\frac{1}{2}\|h\|^2]dt 
            -Y_Th(T, X_T)
        \right]
    \right\}.
\end{align}
The right hand side of \eqref{eq:KS} is the likelihood of the measurement.
The variational form of the smoothing problem seeks a distribution $\tilde\cP$ on the path space that minimizes 
	\begin{align}\label{eq:variational}
		\mE_{\tilde\cP} \left\{\log\frac{d\tilde\cP}{d\cQ^Y}\right\}&=\mE_{\tilde\cP} \left\{\log \frac{d\tilde\cP}{d\cP} -\log\frac{d\cQ^Y}{d\cP}\right\}\\
		&={\rm KL} (\tilde\cP\,\|\,\cP) -\mE_{\tilde\cP} \left\{\log\frac{d\cQ^Y}{d\cP}\right\}. \nonumber
	\end{align}
	
Let $\tilde\cP$ be parametrized by the diffusion process $\tilde X_t$ with dynamics
	\begin{equation}\label{eq:tildeX}
		d\tilde X_t = b(t,\tilde X_t) dt + \sigma(t,\tilde X_t) (u_tdt+dW_t),\quad \tilde X_0 \sim \pi_0.
	\end{equation}
By Girsanov theorem \eqref{eq:girsanovthm}, 
	\begin{equation}\label{eq:KLinit}
		{\rm KL} (\tilde\cP\,\|\,\cP) = \mE \left\{\int_0^T \frac{1}{2}\|u_t\|^2 dt\right\} + {\rm KL} (\pi_0\,\|\,\nu_0).
	\end{equation}
Note that \eqref{eq:KLinit} is slightly different from \eqref{eq:KL} since in the control problem $\cP_u$ and $\cP_0$ share the same initial distribution while \eqref{eq:tildeX} and \eqref{eq:filteringdyn1} don't. 
Plugging \eqref{eq:KS} and \eqref{eq:KLinit} into \eqref{eq:variational} yields
an optimal control-like formulation \cite{KimMeh20}
	\begin{align}\label{eq:controlsmoothing}
 \min_{u,\pi_0} ~&\mE\left\{\int_0^T [\frac{1}{2}\|u_t\|^2+\frac{1}{2\sigma_B^2}\|h(t,\tilde X_t)\|^2] dt \right. \\
 &\left. +\frac{1}{\sigma_B^2}Y_t dh(t,\tilde X_t)-\frac{1}{\sigma_B^2}Y_Th(T,\tilde X_T)\right\}+  {\rm KL} (\pi_0\,\|\,\nu_0) \nonumber
	\end{align}
for the smoothing problem. Apart from an extra term ${\rm KL} (\pi_0\,\|\,\nu_0)$ related to the initial distributions, \eqref{eq:controlsmoothing} coincides with the optimal control problem \eqref{eq:sysSDE}-\eqref{eq:cost}-\eqref{eq:costspecial} if we take 
	\begin{subequations}\label{eq:ginference}
	\begin{align}\label{eq:ginference1}
		g(t,x)dt &= \frac{1}{2\sigma_B^2}\|h(t,\tilde x)\|^2 dt+\frac{1}{\sigma_B^2}Y_t dh(t,\tilde X_t), \\
        \Psi(x) &= -\frac{1}{\sigma_B^2}Y_Th(T,x). \label{eq:ginference2}
	\end{align}
	\end{subequations}
	
\section{Path integral particle smoothing}\label{sec:controlfiltering}
Building on the control formulation of smoothing \eqref{eq:controlsmoothing} and the linear approach to optimal control presented in Section \ref{sec:linearapproach}, we propose a new particle filtering algorithm under the name ``path integral particle filtering (PIPF).'' 

\subsection{Path integral particle smoothing}\label{sec:PIsmoothing}
We begin with the smoothing problem to estimate the posterior distribution $\cQ^Y$ for the system \eqref{eq:filteringdyn} over the time-window $[0,\,T]$. As discussed in Section \ref{sec:duality}, this smoothing problem amounts to an optimal control problem
	\begin{equation}\label{eq:costfiltering}
		\min_{u,\pi_0} \mE\left\{\int_0^T [\frac{1}{2}\|u_t\|^2+g(t,X_t)] dt+\Psi(X_T)\right\}+{\rm KL} (\pi_0\,\|\,\nu_0),
	\end{equation}
where $g$ and $\Psi$ are given in \eqref{eq:ginference}. The only difference to a standard optimal control problem is that the initial distribution $\pi_0$, apart from the control $u$, is also an optimization variable. 

Clearly, this small difference doesn't change the optimal control strategy, which remains to be $u^\star (t,x) = \sigma(t,x)'\nabla \log \varphi(t,x)$ with $\varphi$ as in \eqref{eq:varphipde}. Indeed, the optimal control strategy to an optimal control problem is invariant with respect to the initial condition. Plugging this optimal control into \eqref{eq:costfiltering} we arrive at the optimization over $\pi_0$, which reads
	\begin{equation}
		\min_{\pi_0} \mE_{\pi_0} \left\{-\log\varphi(0,X_0)-\log \nu_0(X_0)+\log\pi_0(X_0)\right\}.
	\end{equation}
Apparently, its optimal solution is
	\begin{equation}
		\pi_0^\star(\cdot) \propto \nu_0(\cdot) \varphi(0,\cdot).
	\end{equation}
Note that $\pi_0^\star$ is exactly the posterior distribution of $X_0$ given the full observation $\{Y_t,\, 0\le t\le T\}$. Thus, to sample from the posterior distribution $\cQ^Y$, one can sample $K$ trajectories $\{X_t^k\}_{k=1}^K$ of the diffusion process \eqref{eq:sysSDE} under optimal control strategy $u^\star$ with initial distribution $\pi_0^\star$. The empirical distribution formed by these $K$ trajectories on the path space $\Omega$ is an estimation of the posterior distribution $\cQ^Y$. Moreover, 
	\begin{equation}
		\frac{1}{K} \sum_{k=1}^K \delta_{X_t^k}
	\end{equation}
is an estimation of the posterior distribution of $X_t$ given the full observation $\{Y_t,\, 0\le t\le T\}$. 

The above sampling strategy requires the exact posterior distribution $\pi_0^\star$ of $X_0$ and the exact optimal control strategy $u^\star$. This can be made possible using path integral control but is still computational demanding. Our strategy to sample from $\cQ^Y$ is to sample trajectories with a suboptimal initial distribution $\pi_0$ and a suboptimal control strategy $u$, and then weight the trajectories through important sampling. More precisely, let $\tilde \cP$ be the measure over the path space $\Omega$ associated with initial distribution $\pi_0$ and a suboptimal control strategy $u$, and $\{X_t^k\}_{k=1}^K$ be $K$ trajectories independently sampled from $\tilde \cP$. By Girsanov theorem, in view of \eqref{eq:KS}, 
	\begin{align*}
		\frac{d\cQ^Y}{d\tilde \cP} &= \frac{d\cQ^Y}{d\cP}\frac{d\cP}{d\tilde\cP}
		\\&\propto \frac{d\nu_0}{d\pi_0} \exp[-S^u(0)]
	\end{align*}
where $S^u$ is defined in \eqref{eq:Su} with $g, \Psi$ as in \eqref{eq:ginference}. Denote the value of $S^u(0)$ along the trajectory $X_t^k$ by $S_k^u(0)$ and define the weights
	\begin{equation}\label{eq:Sw}
		w^k =  \frac{d\nu_0}{d\pi_0}(X_0^k) \exp[-S_k^u(0)].
	\end{equation}
It follows that $\cQ^Y$ can be approximated by the empirical distribution formed by the trajectories  $\{X_t^k\}_{k=1}^K$ and weights $\{w^k\}_{k=1}^K$, that is, 
	\begin{equation}\label{eq:weightappro}
		d\cQ^Y \approx \sum_{k=1}^K \hat w^k \delta_{X_{(\cdot)}^k},
	\end{equation}
where 
	\[
		\hat w^k = \frac{w^k}{\sum_{k=1}^K w^k}
	\]
are the normalized weights. Similarly, the posterior distribution of $X_t$ is approximated by
	\begin{equation}
		\sum_{k=1}^K \hat w^k \delta_{X_{t}^k}.
	\end{equation}

The effectiveness of the above approximation \eqref{eq:weightappro} depends on the variance of the weights $\{w^k\}_{k=1}^K$. This variance reduces to zero when $\pi_0$ and $u$ are optimal, that is, $\pi_0=\pi_0^\star, u = u^\star$. In general, computing the exact optimal solution is too expensive and one has to use a suboptimal solution that is easier to compute. There are many methods that can generate suboptimal controller for \eqref{eq:costfiltering}, including differential dynamic programming (DDP) \cite{JacMay70} and iterative linear quadratic regulator (iLQR) \cite{LiTod04}. One can also start from the original smoothing problem for \eqref{eq:filteringdyn} and adopt suboptimal smoothing methods such as extended Rauch-Rung-Striebel (ERTS) \cite{RauTunStr65}. These suboptimal smoothing methods induce suboptimal $\pi_0$ and $u$ for \eqref{eq:costfiltering}. 

    To summarize, our path integral particle smoothing method consists of a proposal initial distribution $\pi_0$ and a proposal feedback $u$. They should be designed such that the distribution on path space induced by $\pi_0$ and $u$ is an approximation of the posterior distribution $\cQ^Y$. A better proposal implies a better estimation with lower variance. Once the proposal is chosen, we can sample trajectories from the controlled diffusion process \eqref{eq:sysSDE} under the proposal control strategy $u$ with the proposal initial distribution $\pi_0$. The posterior distribution $\cQ^Y$ is then approximated by \eqref{eq:weightappro}.
    
\subsection{Path integral particle filtering}
We next move to the filtering problem. We are interested in the filtering problem of estimating the posterior distribution of $X_t$ conditioning on the past observation $\{Y_\tau, 0\le \tau\le t\}$. More precisely, denote by $\cY_t=\sigma(Y_\tau: 0\le \tau\le t)$ the sigma-field generated by the observation up to time $t$, then the objective of filtering is to estimation $\cP(X_t \in \cdot \mid \cY_t)$.

The path integral particle smoothing algorithm proposed in Section \ref{sec:PIsmoothing} is suitable for smoothing problem over a fixed time window $[0,\,T]$. To use this method for filtering problem where new measurements keep coming in, one naive strategy is to carryout smoothing task over the time window $[0,\, t]$. However, this requires recursively implementing the smoothing algorithm over a larger and larger time window. As $t$ increases, the computational complexity of the smoothing problem grows and will eventually become computationally infeasible. 

We propose to use a sliding window implementation of the smoothing algorithm for filtering. More specifically, consider the smoothing problem over the time window $[t-H, t]$ of size $H>0$. It is equivalent to the optimal control problem
\begin{align}\label{eq:w-costfiltering}
    \min_{u,\pi_{t-H}} \mE\left\{\int_{t-H}^t [\frac{1}{2}\|u_\tau\|^2 +g(\tau, X_\tau)] d\tau+\Psi_t(X_t)\right\}\\ 
    +{\rm KL} (\pi_{t-H}\,\|\,\nu_{t-H}), \nonumber
\end{align}
where $g$ is as in \eqref{eq:ginference1} and 
	\begin{equation}
		\Psi_t(x) = -\frac{1}{\sigma_B^2}Y_th(t,x).
	\end{equation}
The prior distribution $\nu_{t-H}$ for this smoothing problem is the posterior distribution $\cP(X_{t-H}\in \cdot\mid \cY_{t-H})$. Since $\nu_{t-H}$ already accounts for all the observations $\{Y_\tau, 0\le \tau \le t-H\}$, the solution to the smoothing problem \eqref{eq:w-costfiltering} in fact induces the exact posterior distribution over the trajectories $\{X_\tau, t-H\le \tau \le t\}$, conditioned on the full history of observation $\{Y_\tau, 0\le \tau \le t\}$. Thus, by running the smoothing algorithm presented in Section \ref{sec:PIsmoothing} over a fixed-size time window $[t-H, \, t]$, we can obtain the posterior distribution of $\cP(X_t\mid \cY_t)$. 

To implement path integral particle smoothing algorithm over the time window $[t-H, t]$, one needs to evaluate $d\nu_{t-H}/d\pi_{t-H}$ as in \eqref{eq:Sw}. However, in the proposed path integral particle filtering method, the distribution $\nu_{t-H}$ doesn't have a closed-form and is represented by a collection of weighted particles as
	\begin{equation}
		\nu_{t-H} \approx \sum_{k=1}^K w_p^k X_{t-H}^k. 
	\end{equation}
Thus, a natural way to sample trajectories over time interval $[t-H,\, t]$ is to initialize them with $\{X_{t-H}^k\}_{k=1}^K$ and then follow the closed-loop dynamics \eqref{eq:sysSDE} under some suboptimal control policy $u$. With this strategy, the proposal initial distribution $\pi_{t-H}$ satisfies 
	\begin{equation}
		\frac{d\nu_{t-H}}{d\pi_{t-H}}(X_{t-H}^k) \propto w_p^k.
	\end{equation}
	
Let $\{X_\tau^k\}_{k=1}^K$ be the $K$ generated trajectories and $S_k^u(t-H,t)$ be the value of \eqref{eq:Su} of the trajectory $X_\tau^k$ over the time interval $[t-H,\, t]$, then the posterior distribution over the trajectory space conditioning on the past observation $\{Y_\tau, 0\le \tau\le t\}$ is approximated by
	\begin{equation}
		\sum_{k=1}^K \hat w ^k X_{(\cdot)}^k
	\end{equation}
with $\{\hat w ^k\}_{k=1}^K$ being the normalized version of the weights
	\begin{equation} \label{eq:weightH}
		w_p^k \exp [-S_k^u(t-H,t)]. 
	\end{equation}

To see the intuition of \eqref{eq:weightH}, assume that $\nu_{t-H}$ is obtained using the path integral particle smoothing algorithm over the time interval with proposal initial distribution $\nu_0$. Following the arguments in Section \ref{sec:PIsmoothing}, by \eqref{eq:Sw}, we know 
	\[
		\hat w^k \propto \exp[-S_k^u(0, t-H)]
	\]
where $S_k^u(0, t-H)$ is evaluated over some sampled trajectory over the time interval $[0,\,t-H]$. Combining it with \eqref{eq:weightH} we conclude that 
	\[
		\hat w ^k\propto \exp[-S_k^u(0, t-H)]\exp [-S_k^u(t-H,t)] = \exp[-S_k^u(0, t)],
	\]
where $S_k^u(0, t)$ is evaluated over the concatenated trajectory of $X_\tau^k, 0\le \tau\le t-H$ and $X_\tau^k, t-H\le \tau\le t$. Instead of resampling the whole trajectory starting from the very beginning, in the sliding window filtering, all the past weights are recorded in the particle representation of $\nu_{t-H}$ and are combined with the measurement over $[t-H,\, t]$ to estimate the posterior distribution.

\subsection{Algorithm and Implementation}
In this section we provide implementation details of the path integral particle filtering algorithm. Let $0=t_0<t_1<t_2<\cdots$ be a sequence of time discretization points. It can be a constant stepsize discretization, i.e., $t_{j+1}-t_j = \Delta t$, or any other more flexible discretization scheme. Set the sliding window size after time-discretization in the path integral particle filtering algorithm to be $H$ with a slightly abuse of notation. 

The proposed particle filtering algorithm can be divided into two stages. For $t_j, j\le H$, the total number of time steps is less than the window size and thus we use path integral particle smoothing over the time interval $[0, t_j]$ to estimate the posterior distribution $\cP(X_{t_j} |\cY_{t_j})$. When $j>H$, we adopt the path integral particle filtering over the time interval $[t_{j-H}, t_j]$. 

In the sliding window stage over $[t_{j-H}, t_j]$, the choice of particle representation for the prior distribution $\nu_{t_{j-H}}$ ($\cP(X_{t_{j-H}} |\cY_{t_{j-H}})$) is crucial. We use the trajectories generated in the previous step over time interval $[t_{j-H-1}, t_{j-1}]$ and weight them properly to obtain an estimation of $\nu_{t_{j-H}}$. The locations of the particles generated in this way account for the measurement up to time $t_{j-1}$ and are thus match better with the posterior distribution $\cP(X_{t_{j-H}} |\cY_{t_{j}})$, which is the ideal proposal initial distribution \cite{DouBriSen06}. More explicitly, the prior distribution $\nu_{t_{j-H}}$ is updated recursively as follows. Let the particle representation of the prior distribution $\nu_{t_{j-H-1}}$ at the previous step be 
	\begin{equation}\label{eq:priorpre}
		\nu_{t_{j-H-1}} \approx \sum_{k=1}^K w_p^k X_{t_{j-H-1}}^k,
	\end{equation}
and $S_k^u$ be the values of \eqref{eq:Su} evaluated over the sampled trajectories over time window $[t_{j-H-1}, t_{j-1}]$, then 
	\begin{equation}\label{eq:prior}
		\nu_{t_{j-H}} \approx \sum_{k=1}^K w_p^k \exp[-S_k^u(t_{j-H-1},t_{j-H})]X_{t_{j-H}}^k.
	\end{equation}
In the above, to simplify the notation, the normalization for the weight is not displayed explicitly. 

The effective size of the samples decreases much slower than the standard SIR filter. Yet, a resampling step is needed after a long time horizon. For resampling, we start with the samples in \eqref{eq:prior}. We resample them based on the weights
	\begin{equation}
		w_p^k \exp[-S_k^u(t_{j-H-1},t_{j-1})],
	\end{equation}
obtaining new samples $\hat X_{t_{j-H}}^k$. These samples follow approximately the distribution $\cP(X_{t_{j-H}}|\cY_{t_{j-1}})$. With these new samples, the prior distribution $\nu_{t_{j-H}}$ is approximated by
	\begin{equation}
		\nu_{t_{j-H}} \approx \sum_{k=1}^K \exp[S_k^u(t_{j-H},t_{j-1})]\hat X_{t_{j-H}}^k.
	\end{equation}

Once a particle representation of $\nu_{t_{j-H}}=\sum_{k=1}^K w_p^k X_{t_{j-H}}^k$ is arrived, one can start from it and apply path integral particle smoothing over the time window $[t_{j-H}, t_{j}]$. This leads to the particle filtering result
	\begin{equation}
		\cP(X_{t_j} |\cY_{t_j}) \approx \sum_{k=1}^K w_p^k \exp[-S_k^u(t_{j-H},t_{j})]X_{t_{j}}^k.
	\end{equation}

\begin{figure}
	\centering
	\includegraphics[width=0.45\textwidth,height=4.5cm]{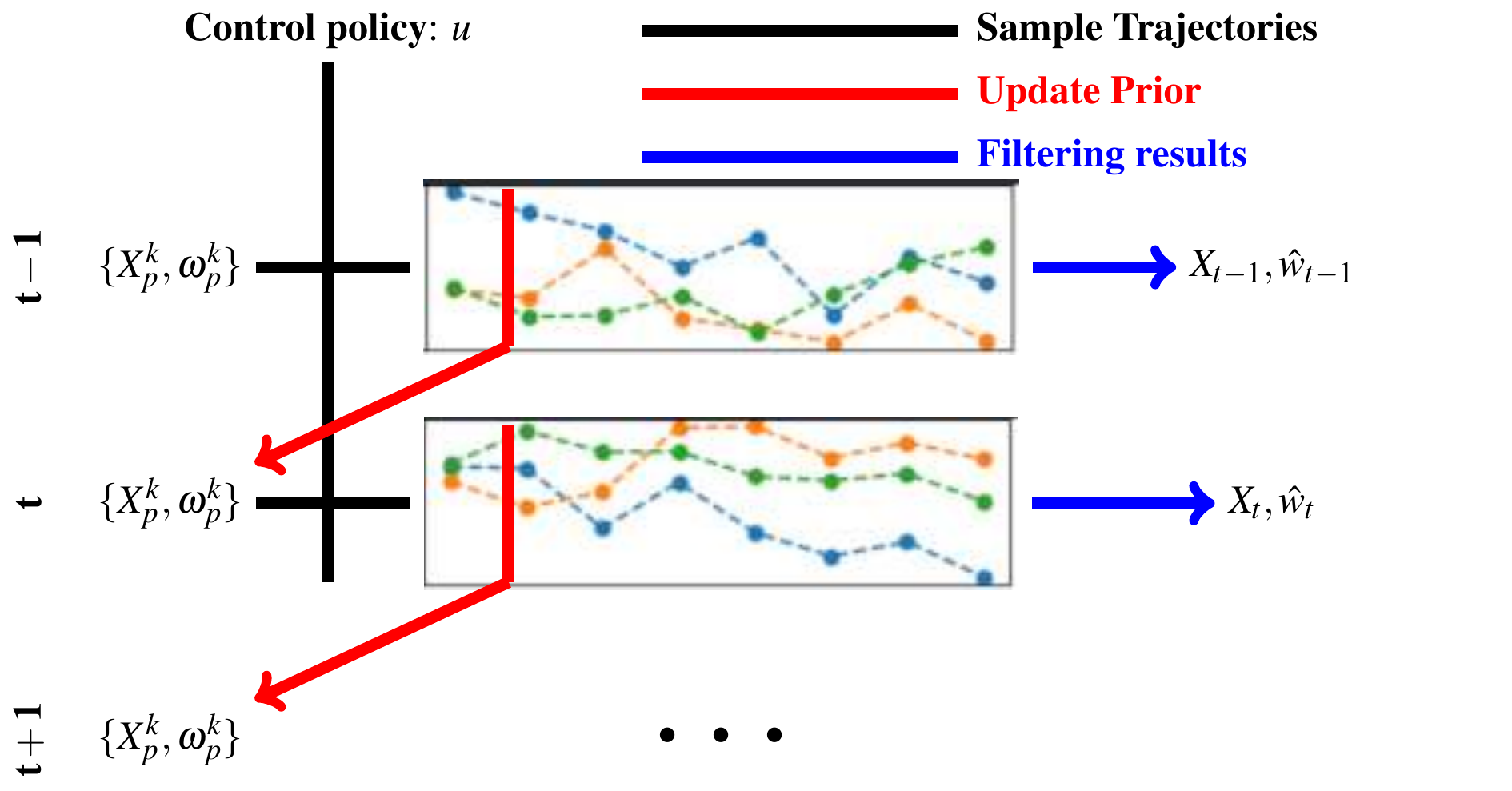}
	\caption{Pipeline diagram of PIPF.}
	\label{fig:pipeline}
\end{figure}

The overall structure of the proposed algorithm is illustrated in the Figure \ref{fig:pipeline}. The full path integral particle filtering algorithm is presented in Algorithm \ref{alg:pi} and a subroutine of it over a given time window is provided in Algorithm \ref{alg:one-step}. The prior distribution at each step is represented by weighted particles $\{X_p^k, w_p^k\}$. The filtering results $\cP(X_{t_j} |\cY_{t_j})$ at the current step is represented by weighted particles $\{X^k, \hat w^k\}$.
\begin{algorithm}
    \caption{Path Integral Particle Filtering (PIPF)}
    \label{alg:pi}
    \begin{algorithmic}
    \STATE \textbf{Input: } $L$: Total Number of time steps 
    \STATE $\{X_p^k\}$: Samples from prior distribution $\nu_0$
    \STATE $\{w_p^k\}$: Weight of samples, all ones
    \STATE $H$: Length of sliding windows
    \FOR{$j \leftarrow 1,\cdots, L$}
    \STATE $j$-th time interval $\leftarrow [\min(0,t_{j-H}), t_j]$
    \STATE $\{X^k\}, \{\hat{w}^k\}, \{X_p^k\}, \{w_p^k\} \leftarrow$ Algorithm~\ref{alg:one-step} for $j$-th time interval
    \ENDFOR
\end{algorithmic}
\end{algorithm}

\begin{algorithm}
    \caption{One Step of PIPF}
    \label{alg:one-step}
    \begin{algorithmic}
        \STATE \textbf{Input: } $b, \sigma, \sigma_B$: System model
        \STATE $g, \Psi$: Cost model as in \eqref{eq:ginference}
         \STATE $\gamma_{thres}$: threshold for resampling
	\STATE $[t_i,\, t_j]$: sliding window
        \STATE $K$: Number of samples
        \STATE $\{X_p^k\}, \{w_p^k\}$: Samples and weight from the previous step
        \STATE \textbf{Output:} $\{X^k\}, \{\hat{w}^k\}$: Filtering results
        \STATE $\{X_p^k\}, \{w_p^k\}$: Samples and weight for the current step
    \FOR{$k \leftarrow 1, \cdots, K$}
        \STATE Sampling $k$-th trajectory $X_\tau^k$ initialized by $X_p^k$.
	\STATE Evaluate the value of $S_k^u$ over the trajectory $X_\tau^k$
	\STATE Filtering samples: $X^k = X_{t_j}^k$
        \STATE Filtering weight: $\hat w^k \propto w^k_p \exp[-S_k^u(t_i,t_j)]$
        \STATE \comment{\# Prior information for usage in next sliding window}
        \STATE Prior sample: $X_p^k = X_{t_{i+1}}^k$
        \STATE Prior weight: $w_p^k \propto w_p^k \exp[-S_k^u(t_i,t_{i+1})]$
    \ENDFOR
    \STATE \comment{\# Check resampling}
    \STATE Effective ratio: $\gamma = \frac{1}{K \sum_{k=1}^K (\hat{w}^k)^2}.$
    \IF{$\gamma < \gamma_{thres}$}
    \STATE $\{X_p^k\} \sim {\rm multinomial} (\{X_p^k\}, \{w_p^k\}\exp[-S_k^u(t_{i},t_j)])$
    \STATE $\{w_p^k\} \propto \exp[S_k^u(t_{i+1},t_j)]$
    \ENDIF
    \end{algorithmic}
\end{algorithm}

The performance of the PIPF algorithm depends on the length $H$ of the sliding window and the choice of proposal suboptimal control $u$.
When $H=1$ and $u\equiv 0$, our algorithm reduces to the standard SIR algorithm as explained further in the following remark.
\begin{rem}
Without any control, i.e. $u=0$,
the  algorithm resembles the SIR particle filter. Indeed, when $u=0$,  the location of the particles $X^k_t$ is only governed by the open-loop dynamics~\eqref{eq:filteringdyn1} similar to the SIR particle filter. And the weights of the particles $w^k_t\propto \exp[-S^0_k(0,t)]$  where 
\begin{align*}
    S^{0}_k(0,t) = &\int_{0}^{t} \frac{1}{2\sigma_B^2}\|h(\tau,X^k_\tau)\|^2d\tau  + \frac{1}{\sigma_B^2}Y_\tau dh(\tau,X^k_\tau) \\&- \frac{1}{\sigma_B^2} Y_t h(t,X^k_t).
\end{align*}
This is precisely the log-likelihood of the observation signal over the time interval $[0,t]$. With a time discretization of the integral $0=t_0<t_1<\ldots<t_J=t$, the weights can be expressed as multiplication of the likelihoods $\prod_{j=1}^J p(Y_{t_j}|X^k_{t_j})$, which is similar to how SIR particle filter updates the weights. 
\end{rem}

\section{Numerical examples}\label{sec:examples}
In this section we present several numerical examples to demonstrate the efficacy of the proposed path integral particle filtering algorithm. In the first example, we test the proposed algorithm in linear filtering problems. In the second example, we consider a nonlinear filtering problem where the optimal filtering can be obtained in closed form and show that PIPF is able to approximate this optimal filtering well. 
\subsection{Linear Filtering exmaples}
\label{ssec:linear-filter-exp}

We first consider the following one-dimensional state space model 
\begin{align*}
    dX_t &= -\kappa X_t dt + dW_t,~X_0 \sim N(m_0, P_0) \\
    dY_t &= X_t dt + \sigma_B dB_t,
\end{align*}
where $\kappa >0$.
 The model corresponds to an \textit{Ornstein-Uhlenbeck} process, whose measurements are corrupted with Gaussian noise. The posterior distribution is Gaussian, that is, $p(X_t | \cY_t)=N(X_t| m_t, P_t)$ with
\begin{align*}
    dm_t &= -\kappa m_t dt + P_t(dY_t - m_t dt),\\
    \frac{dP_t}{dt} &= -2\kappa P_t - P_t^2 + 1.
\end{align*}
Three filtering algorithms are compared for this problem: (i) the sequential importance resampling~(SIR) particle filter~\cite{DouBriSen06}; (ii) the path integral particle filter with zero controller~(PIPF-zero) $u\equiv 0$; (iii) the path integral particle filter with linear quadratic regulator controller~(PIPF-LQR). LQR is designed based on the cost function given in \eqref{eq:ginference}. 
The simulations are executed for $L=600$ times-steps with step-size $\Delta t = 0.01$. Both employ a sliding window size $H=20$. All three algorithms use $K=500$ particles. 

\begin{figure}[h]
	\centering
	\includegraphics[width=0.23\textwidth,height=4cm]{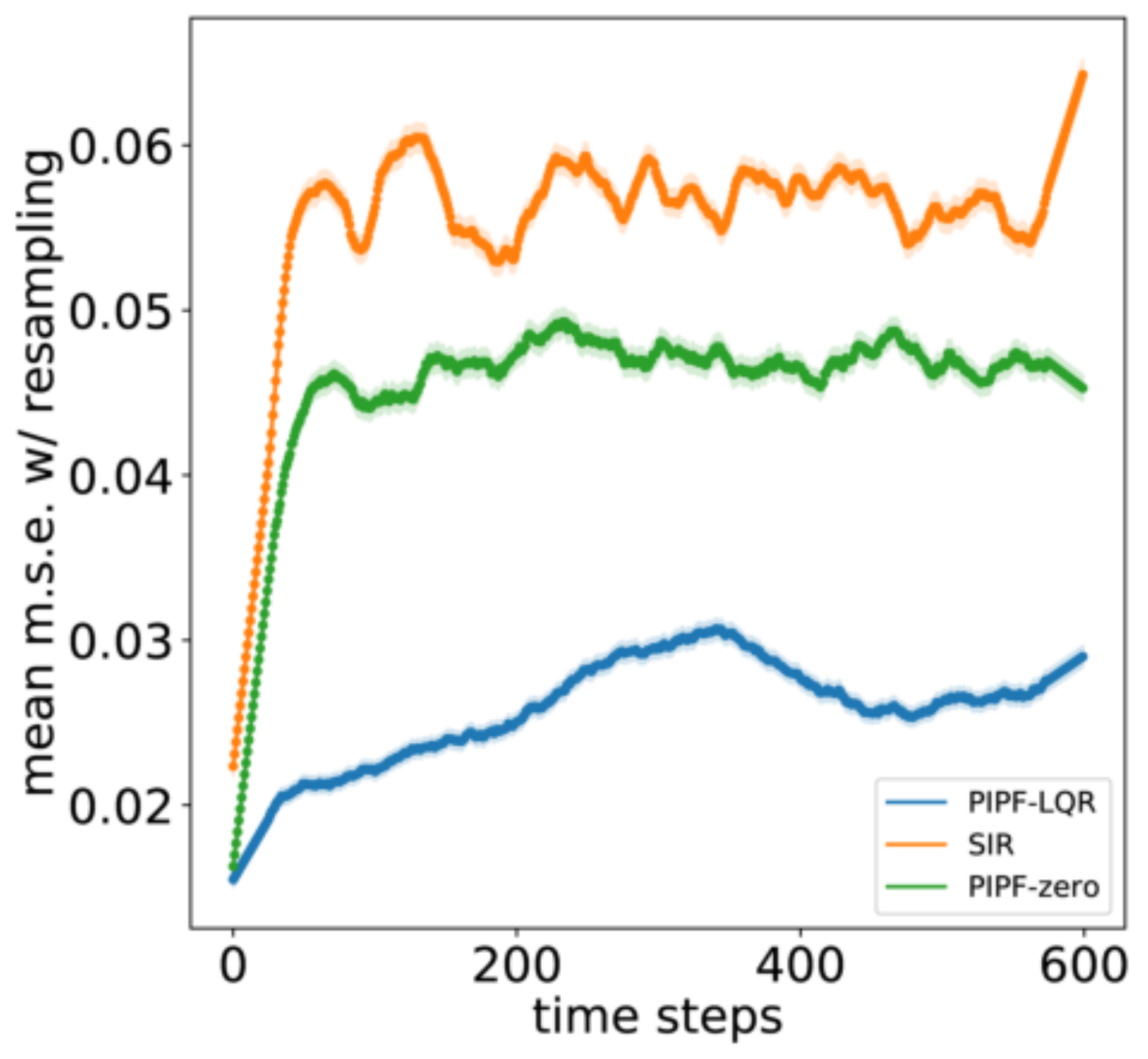}
	\includegraphics[width=0.23\textwidth,height=4cm]{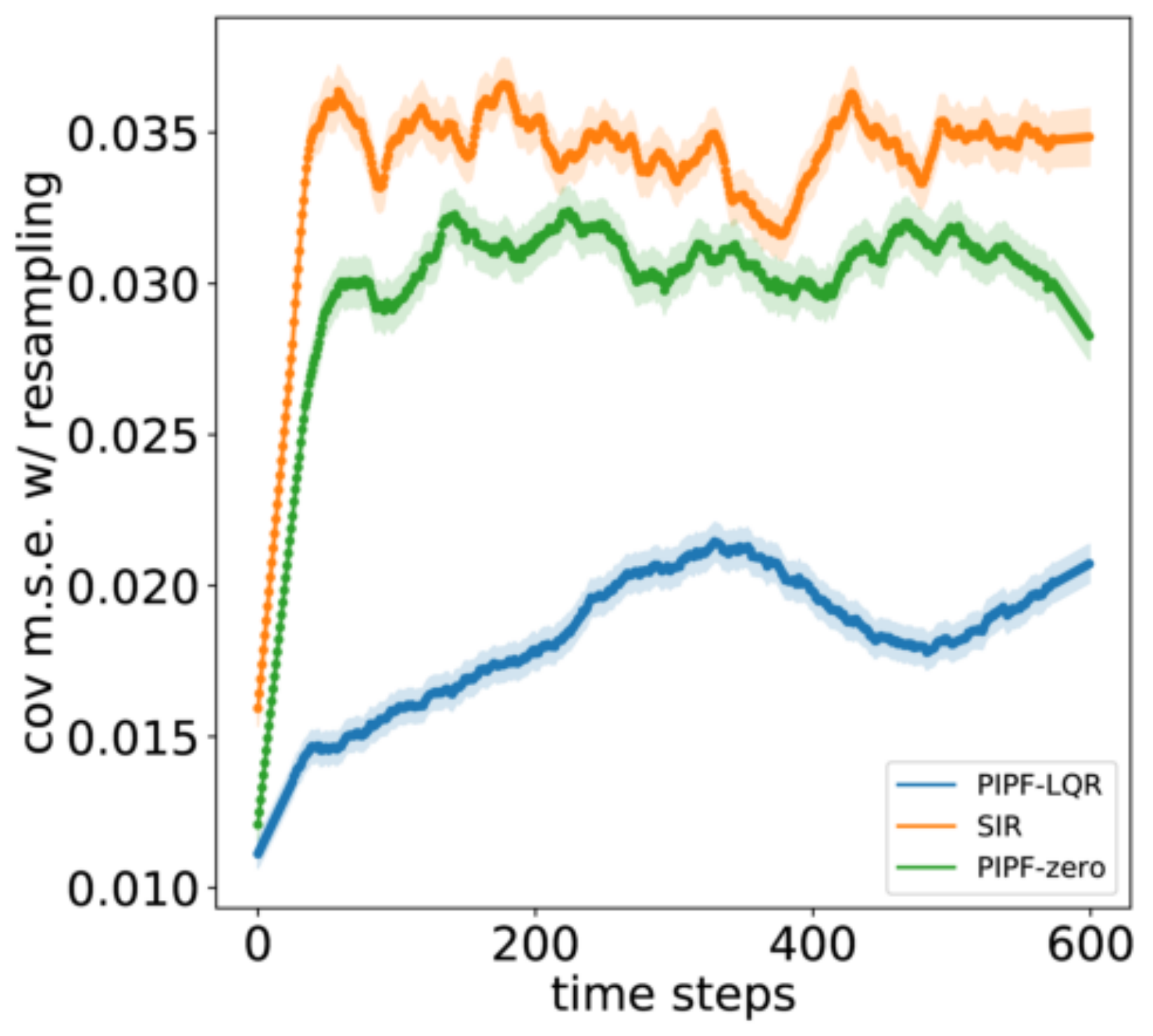}
	\caption{Performance comparison w/ resampling}
	\label{fig:ou}
\end{figure}

\begin{figure}[h]
	\centering
	\includegraphics[width=0.23\textwidth,height=4cm]{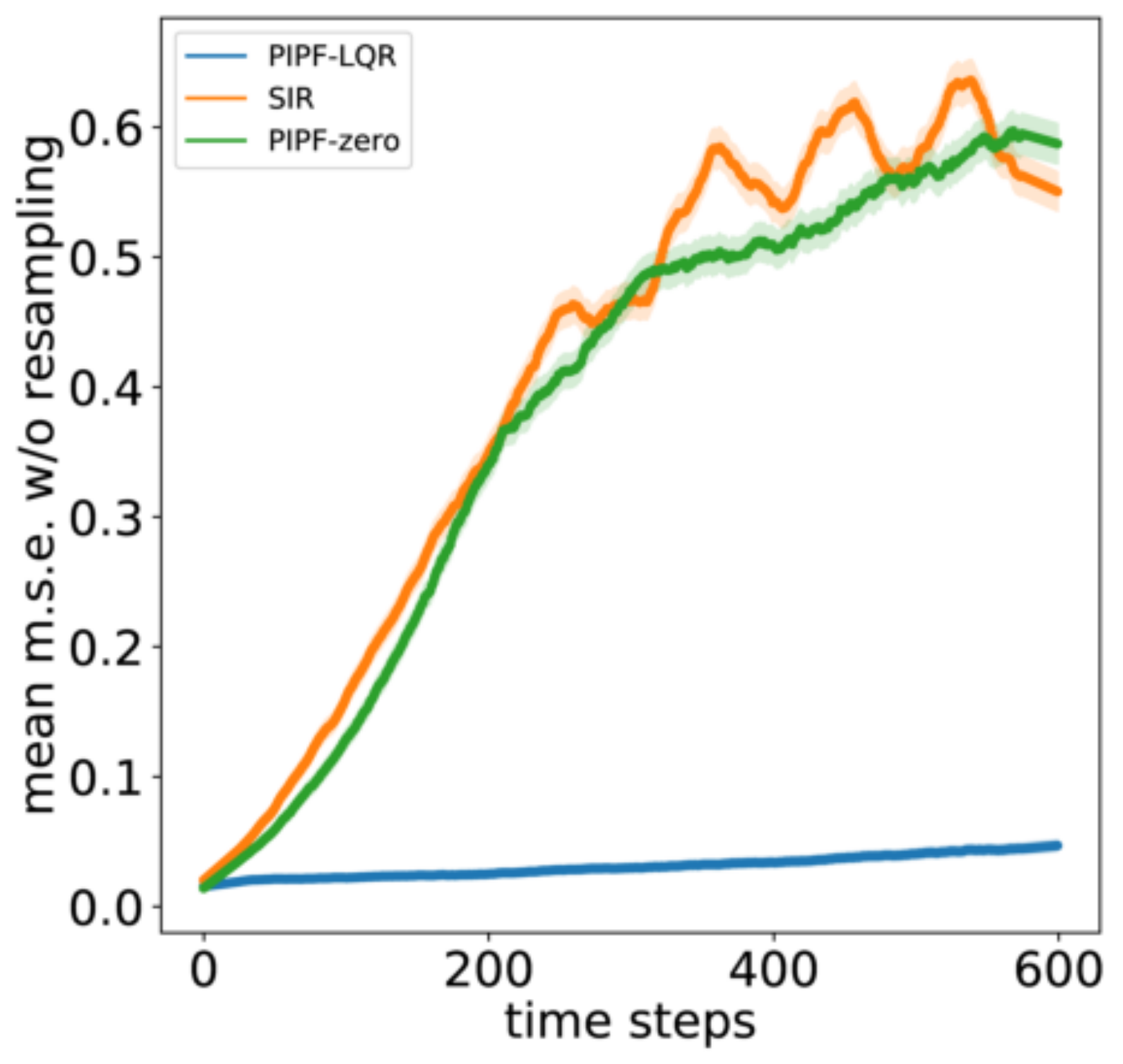}
	\includegraphics[width=0.23\textwidth,height=4cm]{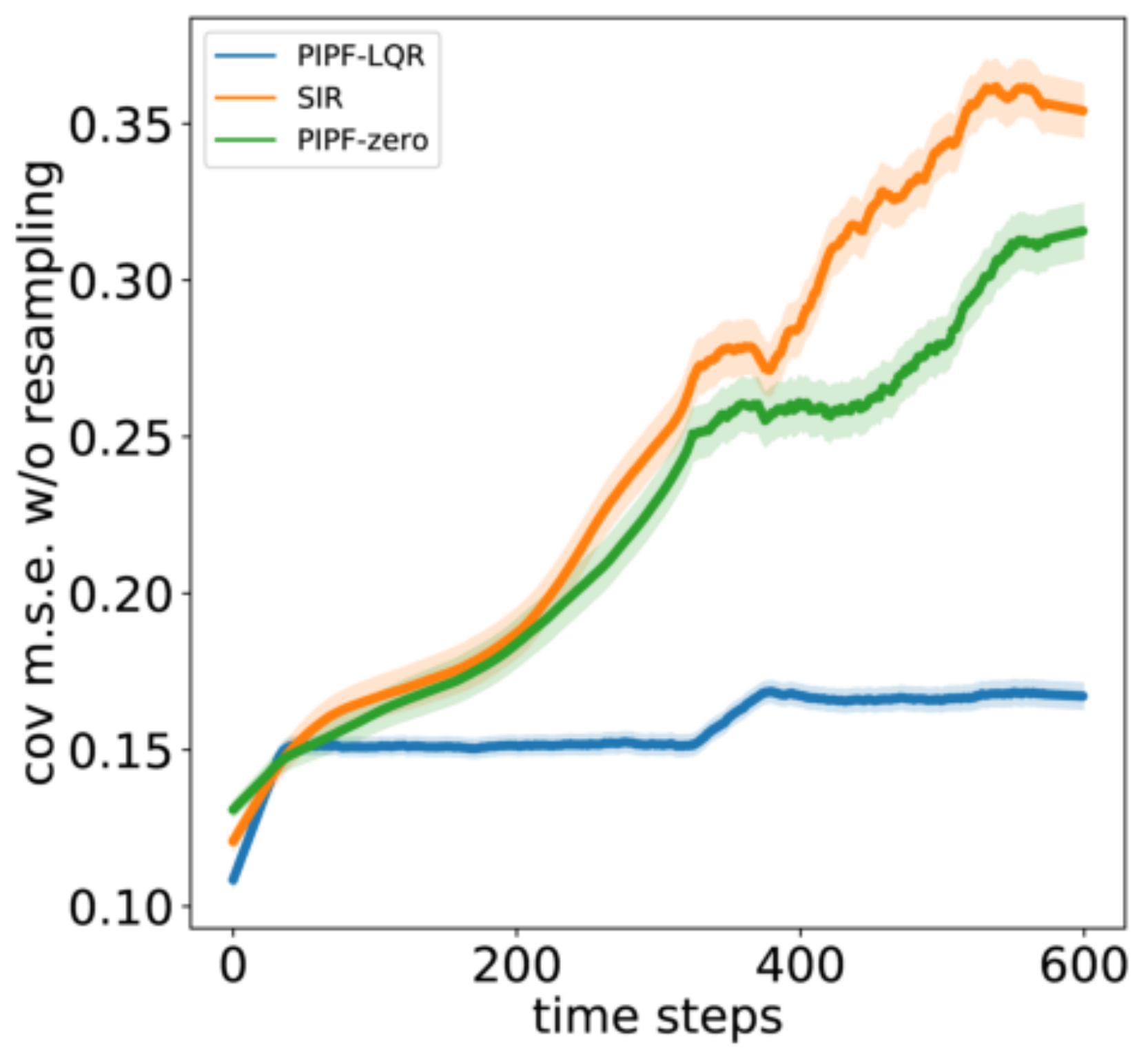}
	\caption{Performance comparison w/o resampling}
	\label{fig:ouwo}
\end{figure}

\begin{figure}[h]
	\centering
	\includegraphics[width=0.34\textwidth,height=4cm]{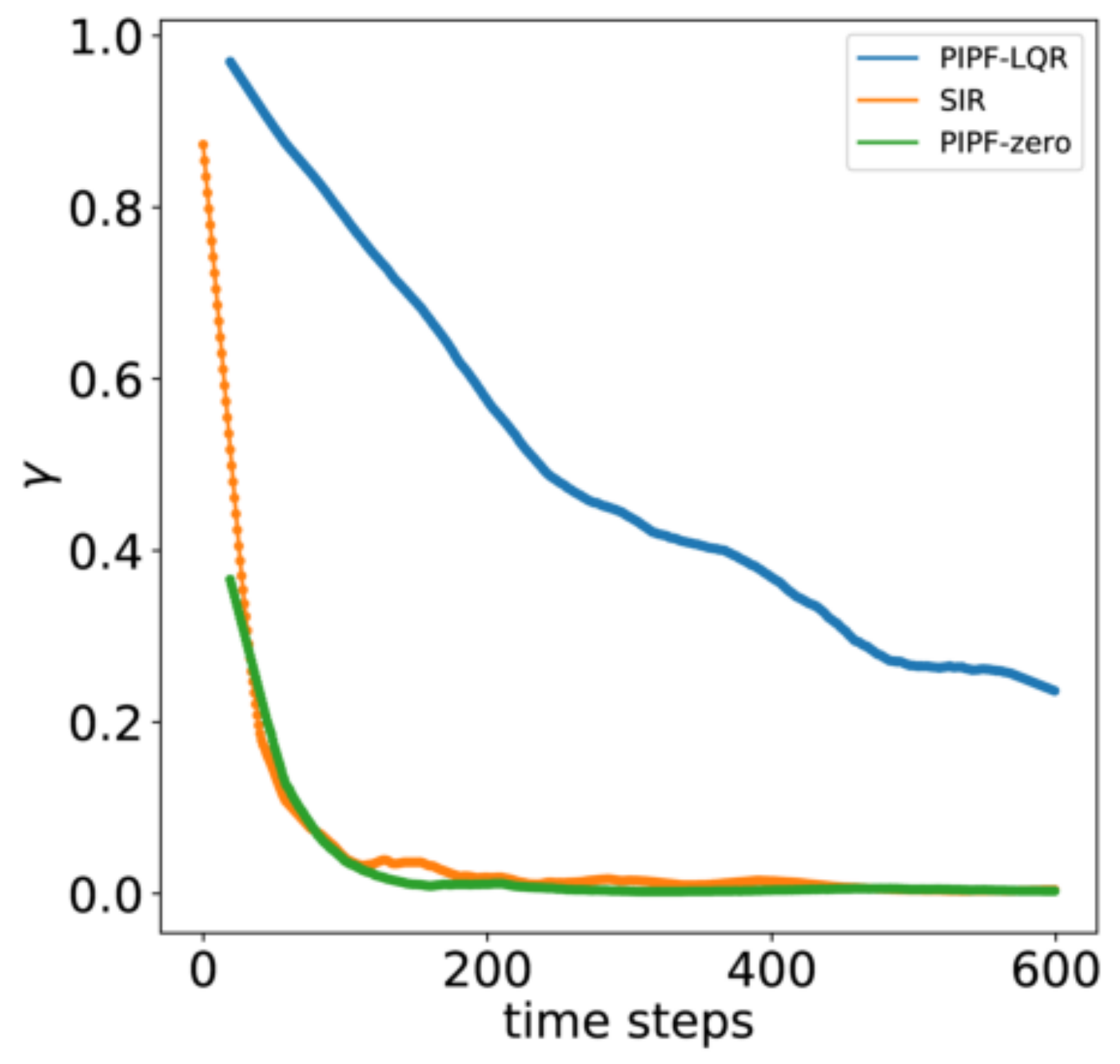}
	\caption{Comparison of effective ratio}
	\label{fig:ouratio}
\end{figure}

It is well known that degenerated particles will lead a low effective ratio and resampling procedure can help ease the notorious degeneracy problem~\cite{DouBriSen06}. 
We benchmark performance in term of the mean squared errors~(m.s.e) between closed-form mean (covariance) and estimated mean (covariance) with resampling and without resampling. The results are depicted in Figure~\ref{fig:ou} and Figure \ref{fig:ouwo} respectively. Another quantity to measure the quality of the particles is the effective ratio $\gamma$. This is depicted in Figure \ref{fig:ouratio} in the absence of resampling. In the experiments, each algorithm is repeated for 50 trials with different random seeds. The solid curves represent the mean of 50 trials and the shaded regions represent the corresponding $1\times$ standard deviation.
From these results we can clearly see the advantage of the path integral particle filtering. Compared with SIR, the introduction of the sliding window in PIPF-zero, which only requires negligible memory footprint, can help reduce the bias of the estimation. The significant improvement of PIPF-LQR owes to the high effective ratio and high-quality particles controlled by the optimal control policy.

\begin{figure}[h]
	\centering
	\includegraphics[width=0.23\textwidth,height=4cm]{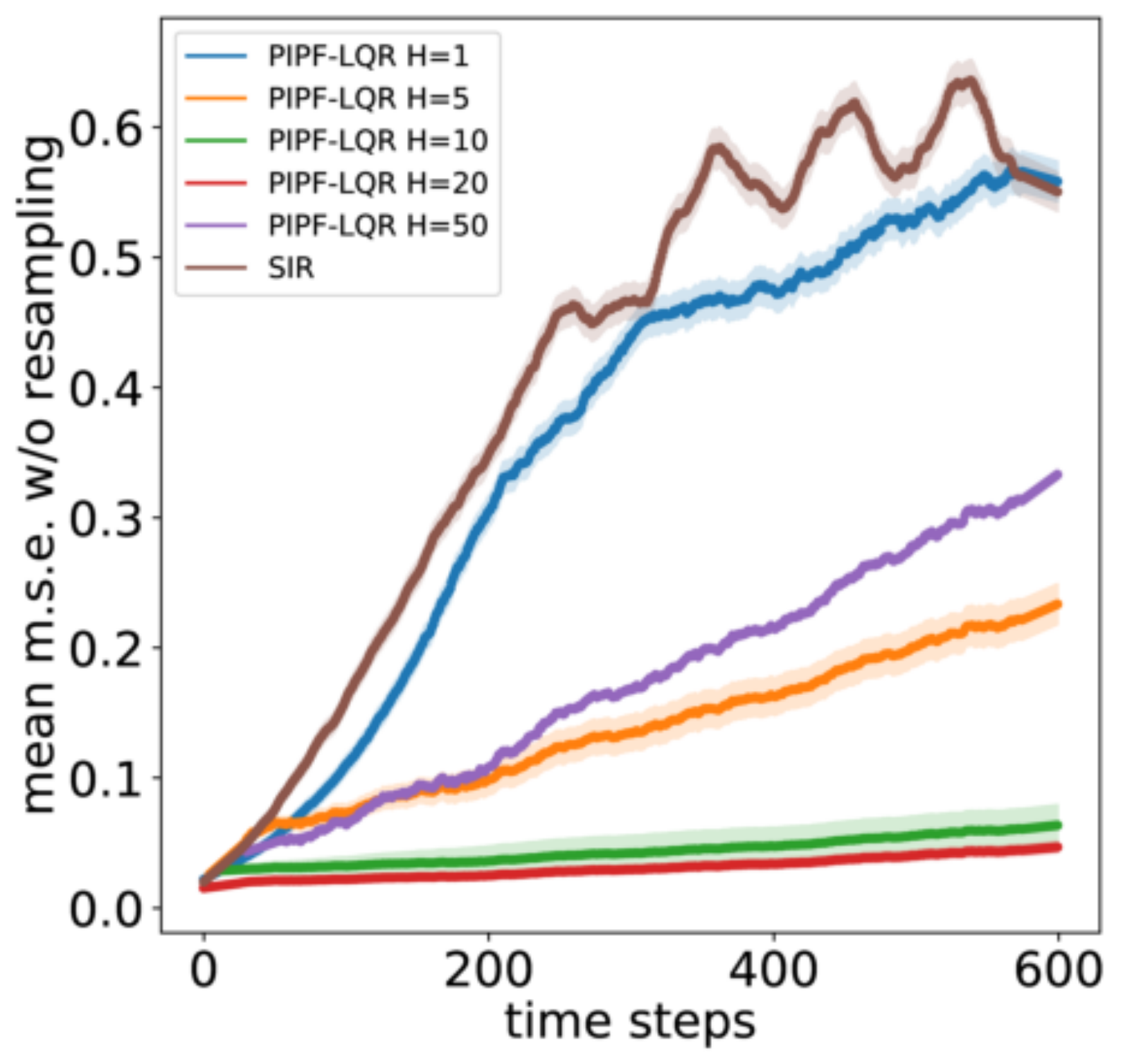}
	\includegraphics[width=0.23\textwidth,height=4cm]{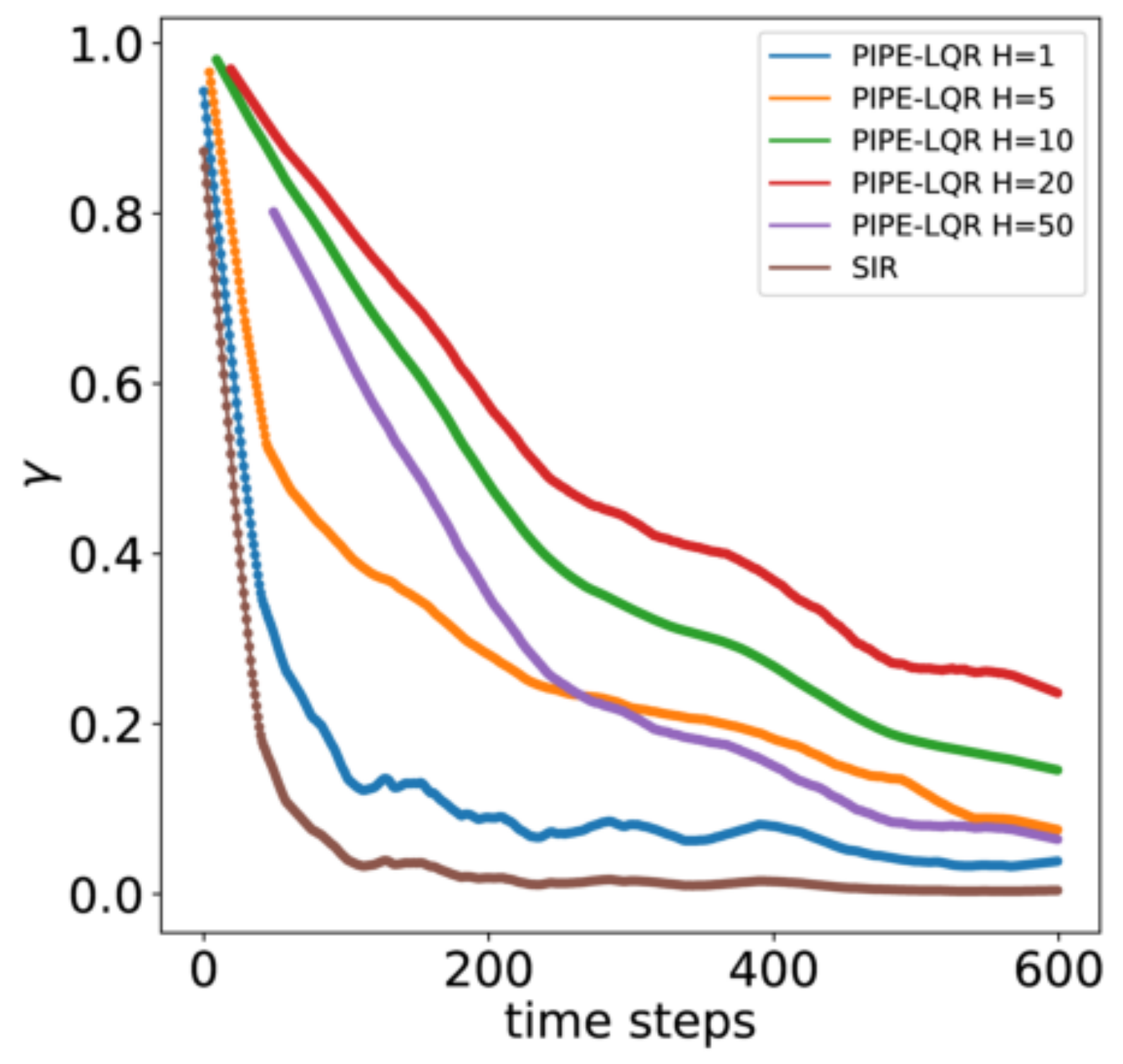}
	\caption{Performance comparison for different $H$.}
	\label{fig:ouH}
\end{figure}

To investigate the influence of the size of the sliding window, we test the PIPF algorithm (with LQR proposal) with different $H$ value. The performance in terms of m.s.e. (of means) with resampling and effective ratio without resampling is depicted in Figure \ref{fig:ouH}. We notice that increasing $H$ improves the filter performance while very large $H$ may deteriorate the results. Parameters other than $H$ are kept the same in this comparison.

\begin{figure}
    \centering
    \includegraphics[width=0.23\textwidth,height=4cm]{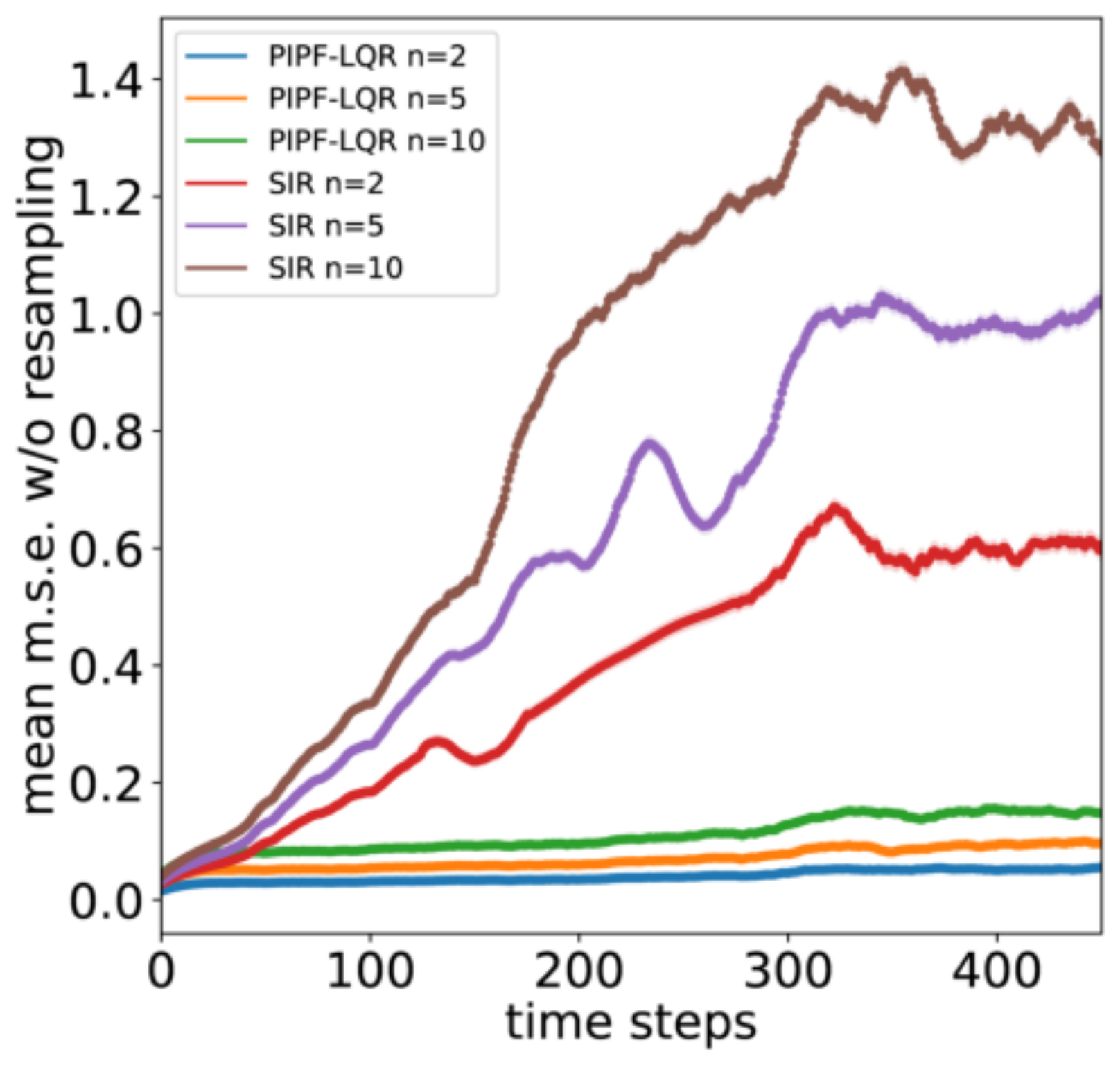}
	\includegraphics[width=0.23\textwidth,height=4cm]{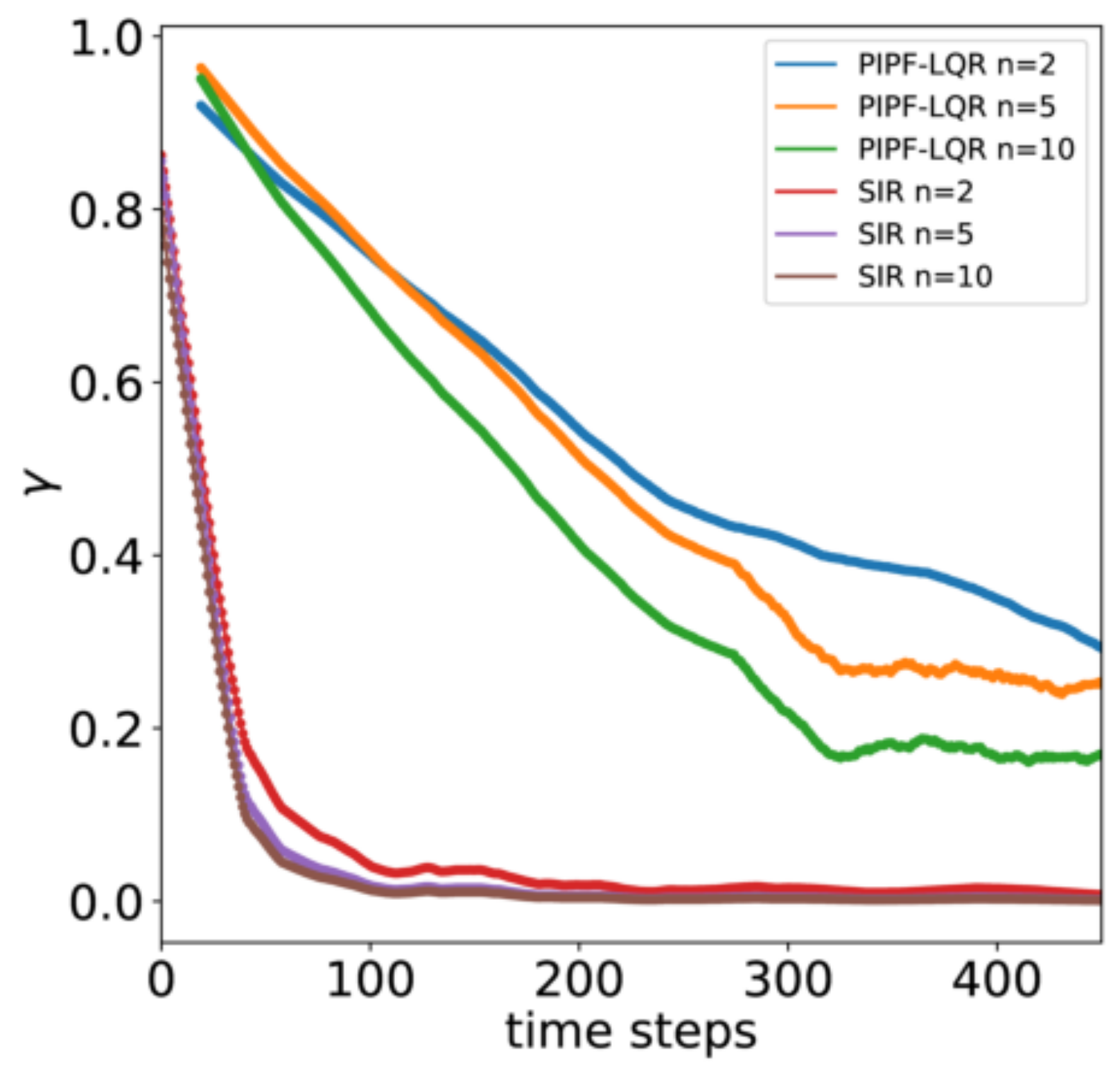}
    \caption{Performance comparison for different $n$.}
    \label{fig:oudim}
\end{figure}

To study the effects of problem dimension on the filtering performance, we consider the multi-dimensional state space model
\begin{align*}
    dX_t &= A X_t dt + dW_t, ~ X_0 \sim N(m_0, P_0) \\
    dY_t &= C X_t dt + \sigma_B dB_t,
\end{align*}
where $A\in\mR^{n\times n}, C\in \mR^{p\times n}$ represent the dynamics matrix and output matrix respectively. In our experiments, $A, C$ are randomly generated for each dimension size $n$ and are then fixed. In Figure \ref{fig:oudim} we display the performance of PIPF-LQR and SIR for different $n$ in terms of m.s.e (of means) and effective ratio without resampling. Clearly, PIPF-LQR scales better than SIR as $n$ increases.

\subsection{Nonlinear Filtering exmaple}
\label{ssec:nonlinear-filter-exp}
We next evaluate our PIPF algorithm on the Benes filter problem \cite{TagMehMey20}
\begin{align*}
    dX_t &= \mu \sigma_W \tanh(\frac{\mu}{\sigma_W} X_t)dt + \sigma_W dW_t, \quad X_0 = x_0, \\
    dY_t &= (h_1X_t + h_1 h_2)dt + dB_t, 
\end{align*}
where $x_0$ is a given constant, and the constants $\mu, \sigma_W, h_1, h_2 \in \mathbb{R}$ are parameters. Regardless being a nonlinear filtering problem, its posterior distribution has an analytical expression
\begin{equation}
    \omega_t N(a_t-b_t, \sigma_t^2) 
    + (1-\omega_t) N(a_t + b_t, \sigma_t^2), 
\end{equation}
where
\begin{align*}
    a_t        & =  \sigma_W \Psi_t \tanh(h_1 \sigma_W t)+\frac{h_2+x_0}{\cosh(h_1) \sigma_W t} -h_2, \\
    b_t        & =  \frac{\mu}{h_1} \tanh(h_1 \sigma_W t), \quad
    \sigma_t^2 =  \frac{\sigma_W}{h_1} \tanh(h_1 \sigma_W t) ,                                      \\
    \Psi_t^2   & = \int_0^t \frac{\sinh(h_1 \sigma_W s)}{\sinh(h_1 \sigma_W t)} dY_s ,\quad \omega_t =  \frac{1}{1+e^{\frac{2a_tb_t}{\sigma_W} \coth(h_1 \sigma_W t)}}.
\end{align*}

In the experiments, we use the model parameters $\mu=1, h_1=1, h_2 = 0, \sigma_W=1, x_0=-5.0.$ The simulation is carried out for $L=6000$ time steps with step-size $\Delta t = 0.001$. 
We compare the performance of our PIPF algorithm and SIR. Since the dynamics is nonlinear, we use a suboptimal control policy, iterative linear quadratic regulator (iLQR) for PIPF. iLQR \cite{LiTod04} approximates the nonlinear dynamics by linearizing it around a nominal trajectory and the cost by a quadratic function, yielding a LQR problem. The iLQR algorithm then solves the resulting LQR problem.
The nominal trajectory is calculated by minimizing control cost Equation~\eqref{eq:costfiltering} locally under noise-free version of the dynamics \eqref{eq:sysSDE}. 
In PIPF-iLQR, the linearization for iLQR is done at the beginning of each sliding window. 

\begin{figure}[h]
    \centering
    \includegraphics[width=0.45\textwidth,height=4cm]{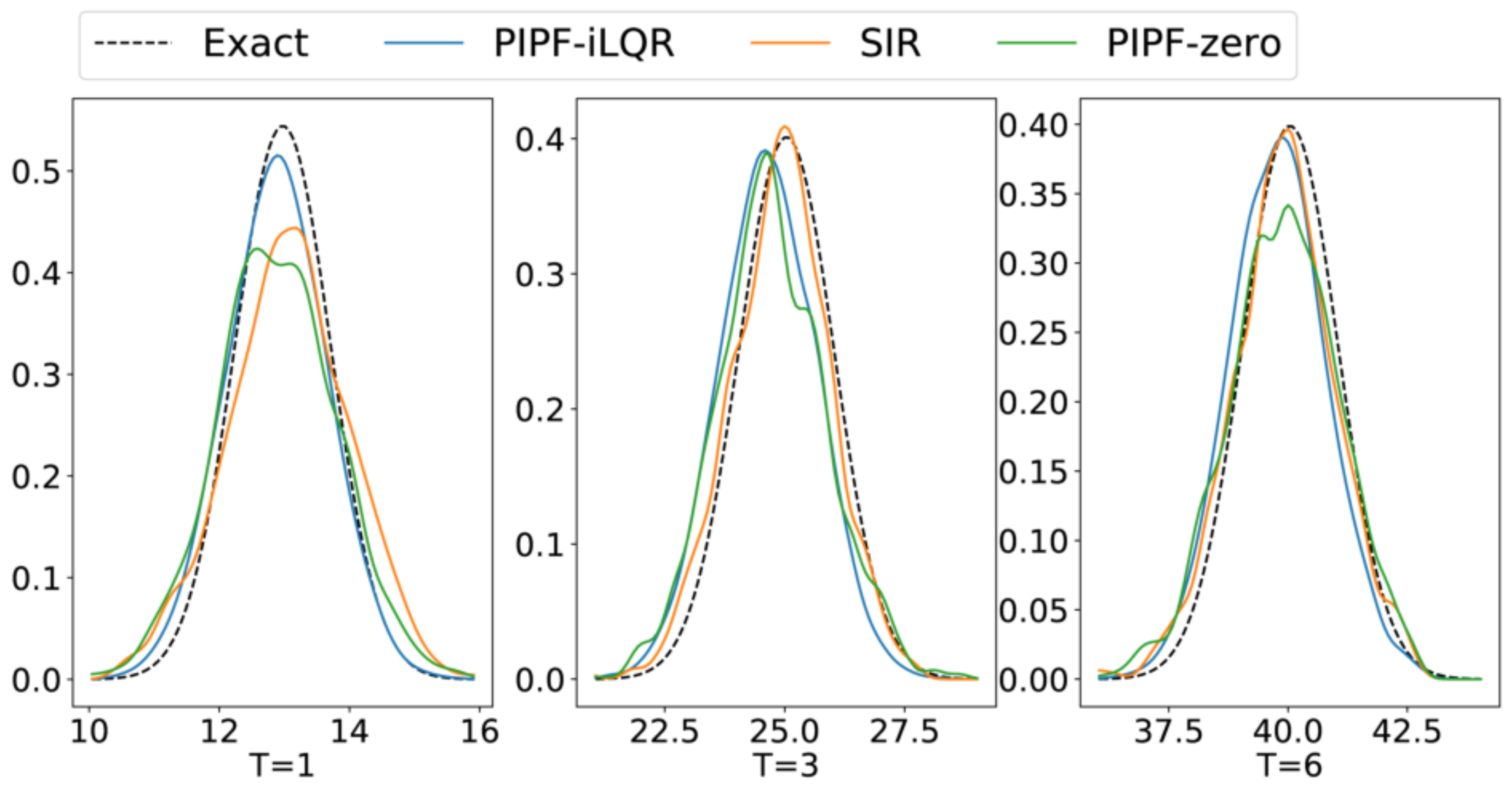}
    \caption{Estimated posterior distribution vs ground truth.}
    \label{fig:benes}
\end{figure}

\begin{figure}[h]
    \centering
    \includegraphics[width=0.23\textwidth,height=4cm]{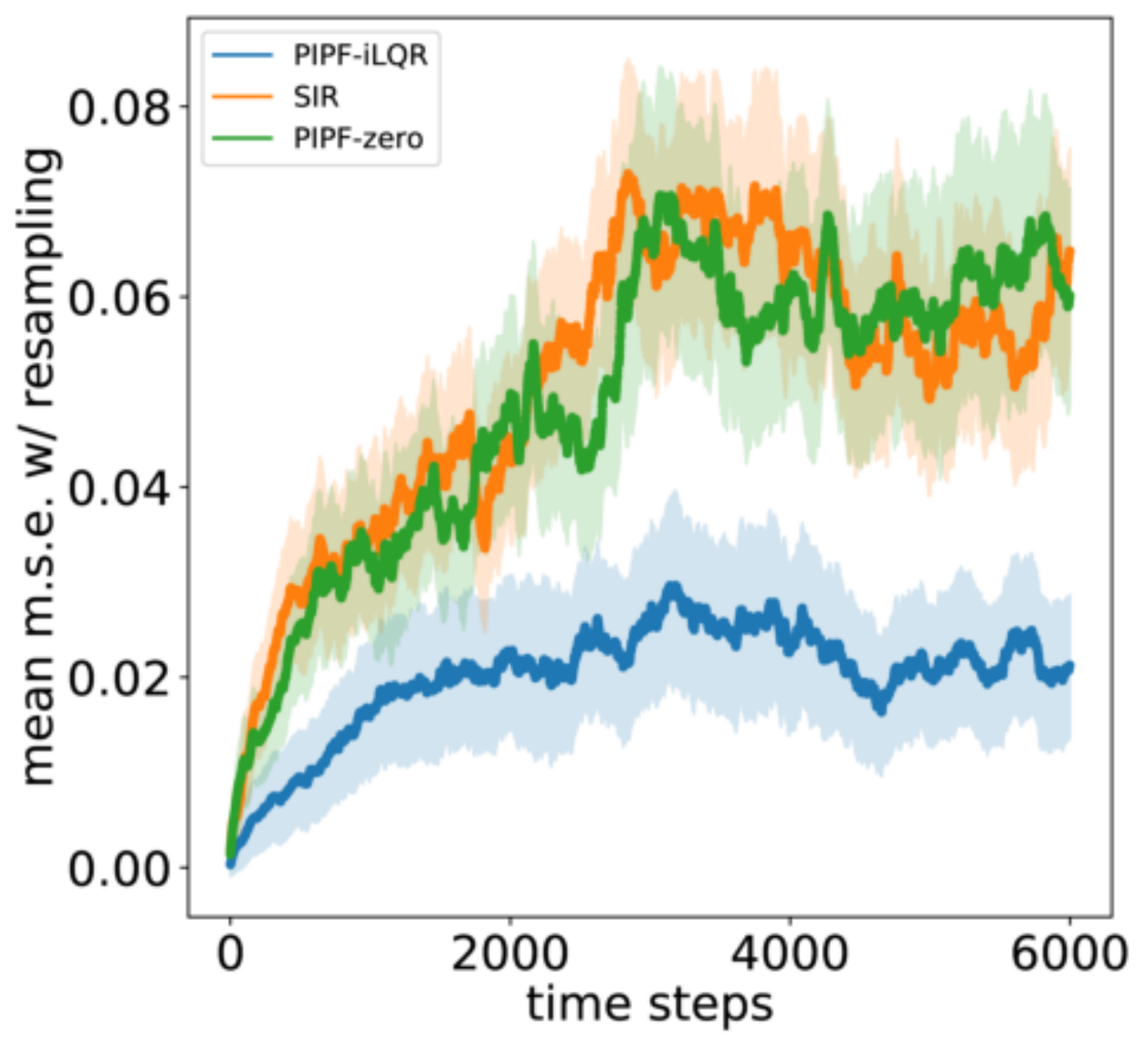}
    \includegraphics[width=0.23\textwidth,height=4cm]{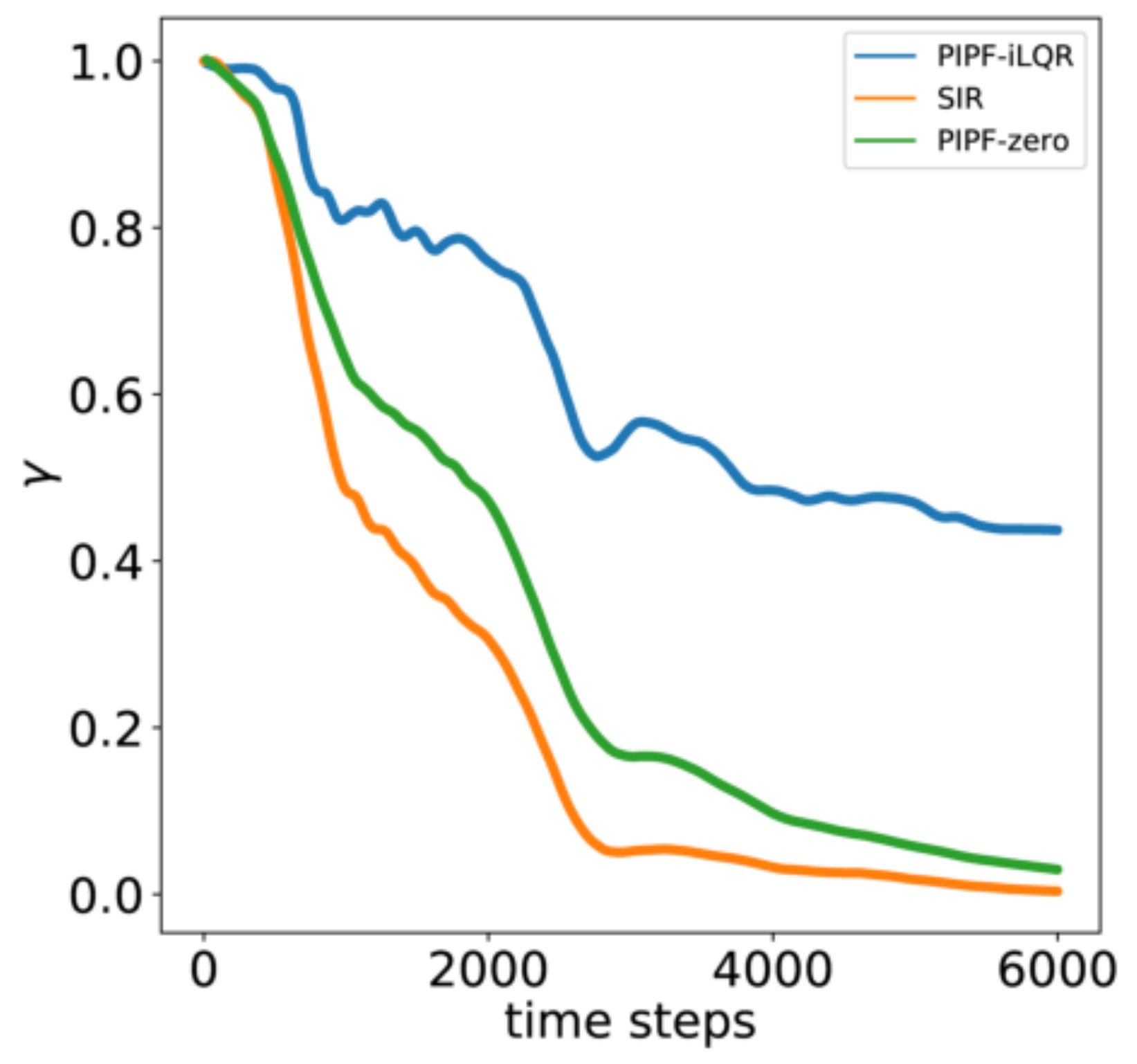}
    \caption{Performance comparison between PIPF and SIR. }
    \label{fig:benesmse}
\end{figure}

The results are displayed in Figure~\ref{fig:benes}-\ref{fig:benesmse}. The sliding window size for PIPF-iLQR and PIPF-zero are set to be  $H=10$. In particular, Figure \ref{fig:benes} displays the estimated posterior distributions of the state at several time points. These distributions are approximated using KDE density estimator \cite{DavLiiPol11} with bandwidth $0.2$. In Figure \ref{fig:benesmse}, we show the m.s.e of the means with resampling and the effective ratio without resampling. 
From the experiments we see that even though the optimality of controller is not promised, PIPF-iLQR still outperforms other algorithms. 

\section{Conclusion}\label{sec:conclusion}

In this paper, buliding on the duality between optimal estimation and optimal control theory, we developed a novel particle filtering algorithm. This algorithm has several distinguish features compared with standard particle filtering algorithms, including high effective ratio and the ability to update samples in the past so as to improve robustness. Our algorithm can also be combined with most existing filtering algorithms such as EKF and UKF to improve their performance. In the future, we plan to extend our algorithm to tackle filtering for diffusion processes with jumps.

{
\bibliographystyle{IEEEtran}
\bibliography{./main}
}
\end{document}